\documentclass[12pt,a4paper]{amsart}
\usepackage{times,amssymb}

\newtheorem{thm}{Theorem}
\newtheorem{lemma}{Lemma}
\newtheorem{proposition}{Proposition}
\newtheorem*{remark}{Remark}

\newcommand{\R}{{\mathbb R}}

\numberwithin{equation}{section}

\begin{document}

\title[Unavoidable configurations of
balls]{Configurations of balls in Euclidean space that Brownian
motion cannot avoid}

\author{Tom Carroll}
\address{Department of Mathematics\newline\indent
University College Cork, Ireland}
\email{t.carroll@ucc.ie}

\author{Joaquim Ortega-Cerd\`a}
\address{Departament de Matem\` atica Applicada i An\`alisi
\newline \indent
Universitat de Barcelona, Gran Via 585, 08007-Barcelona, Spain}
\email{jortega@ub.edu}

\thanks{The second author is  supported by projects MTM2005-08984-C02-02 and
2005SGR00611}

\date{Final version of 13 January 2006 as submitted to Annales Acad.\ Sci.\ Fennic{\ae}}
\maketitle
\begin{abstract} \noindent We consider a collection of
balls in Euclidean space and the problem of determining if Brownian
motion has a positive probability of avoiding all the balls
indefinitely.
\end{abstract}

\section{Introduction}
\noindent We write $\overline{B}(c,r)$ for the closed ball in $\R^d$
with centre $c$ and radius $r$, and write $S(c,r)$ for the sphere of
that centre and radius. We consider a region $\Omega$ that is formed
by removing a countable collection of \textsl{non-overlapping\/}
closed balls from $\R^d$. Thus
\[
\Omega = \R^d \setminus \bigcup_{n=1}^\infty
\overline{B}(c_n,r_n),
\]
and we assume, for convenience, that $0$ lies in $\Omega$. We say
that such a collection of balls is \emph{avoidable} if there is a
positive probability that Brownian motion in $\R^d$, starting from
0, never hits any of the balls. Thus the collection of balls is
avoidable if the balls do not have full harmonic measure w.r.t.\ the
domain $\Omega$, or if infinity has positive harmonic measure
w.r.t.\ $\Omega$. We address the problem of obtaining a geometric
characterization of avoidable configurations of balls.

The genesis of this problem is to be found in the paper of
Ortega-Cerd\`a and Seip \cite{OCS}. Motivated by a question of
Akeroyd \cite{AKER}, the analogous problem in the setting of the
unit disk was completely solved when the centres of the disks that
are removed form a uniformly dense sequence.

In the plane it is possible to hide infinity from $0$ with a single
disk. This reflects the fact that Brownian motion in the plane is
recurrent and that the sphere $S(c,r)$ has full harmonic measure
with respect to $\R^2 \setminus \overline{B}(c,r)$. For this reason,
our results are set in Euclidean space of dimension three or more,
in which Brownian motion is transient. It is helpful to bear in mind
that, in dimension three or more, $( r/\vert c\vert)^{d-2}$ is the
harmonic measure at $0$ of the sphere $S(c,r)$ with respect to the
domain $\Omega = \R^d\setminus \overline{B}(c,r)$. In fact, the
harmonic measure of this sphere at $x$ is $u(x) = ( r/\vert x-c\vert
)^{d-2}$.

\begin{proposition}\label{t1}
We suppose that $d\geq 3$. If
\begin{equation}
\sum_{n=1}^\infty \left( \frac{r_n}{\vert c_n\vert} \right)^{d-2} <
\infty
\label{0}
\end{equation}
then the collection of balls $\{\overline{B}(c_n,r_n)\}_{n\geq 1}$
is avoidable.
\end{proposition}

In order to avoid situations in which a number of small balls packed
very close together can contribute significantly to the sum in
\eqref{0} but contribute relatively little to the overall harmonic
measure, we now require a separation condition on the balls:
\begin{equation}\label{separation}\tag{S}
\text{there is a
positive number $\epsilon$ such that $\vert c_n-c_m\vert \geq
2\epsilon$ for $n\neq m$. }
\end{equation}
\begin{thm}\label{t2}
We suppose that $d\geq 3$. We assume the separation condition
\eqref{separation}
and that there is a number $M$ such that
\begin{equation}
r_n^{d-2}\vert c_n\vert^{2} \leq M \ \text{for }\ n \geq 1.
\label{1}
\end{equation}
If the collection of balls $\{\overline{B}(c_n,r_n)\}_{n\geq 1}$ is
avoidable then
\begin{equation}
\sum_{n=1}^\infty \left( \frac{r_n}{\vert c_n\vert} \right)^{d-2} <
 \infty.
\label{2}
\end{equation}
\end{thm}

The solid angle subtended by the sphere $S(c,r)$ at $0$ is
proportional to $(r/\vert c \vert)^{d-1}$. The appropriate version
of Akeroyd's question in the present setting is whether there is an
avoidable sequence of balls for which  the sum $\sum_n (r_n/\vert
c_n \vert)^{d-1}$ is finite.  If so, it is possible to hide infinity
from the origin from the point of view of harmonic measure even
though geometrically there is a clear line of sight to infinity
except for a set of directions on the sphere $S^d$ of arbitrarily
small $(d-1)$-dimensional measure. Consider $m^{d-1}$ balls of
radius $\rho_m$, with $\rho_m^{d-2} = 1/m^2$, arranged evenly on the
sphere $S(0,m)$, this for each integer $m$ greater than some large
$m_0$. These balls will be non-intersecting and separated, and
\eqref{1} will hold since $\rho_m^{d-2}m^2=1$. But \eqref{2} does
not hold: in fact
\[
\sum_n\left(\frac{r_n}{\vert c_n\vert}\right)^{d-2} =
\sum_{m=m_0}^\infty m^{d-1}\left(\frac{\rho_m^{d-2}}{m^{d-2}}\right)
= \sum_{m=m_0}^\infty \frac{1}{m}.
\]
By Theorem~\ref{t2}, the collection of balls is unavoidable. Even
so, $\sum_n (r_n/\vert c_n \vert)^{d-1}$ is finite.

We will now consider a more regular configuration of balls. We say
that the balls are \emph{regularly located} if (i) the separation
condition \eqref{separation} is satisfied, (ii) the balls are uniformly dense,
in
that there is a positive $R$ such that any ball $B(x,R)$ contains at
least one centre $c_n$, (iii) the radius of any ball depends only on
the distance from the ball's centre to the origin, with $r_n =
\phi(\vert c_n\vert )$ where $\phi$ is a decreasing positive
function.

\begin{thm}\label{t3}We suppose that $d\geq 3$ and that the balls
$B(c_n,r_n)$, $n \geq 1$, are regularly located. Then the collection
of balls is avoidable if and only if
\[
\int^\infty r \phi(r)^{d-2}\,dr < \infty.
\]
\end{thm}

Theorem~\ref{t2} is a partial converse to Proposition~\ref{t1} in
that if the radii of the balls decrease sufficiently rapidly then
the collection of balls is avoidable only if \eqref{0} holds.
Theorem~\ref{t3} will be proved by showing that condition \eqref{1}
is automatically satisfied if the collection of balls is both
regularly located and avoidable. Hence these results do not give
rise to a configuration of separated balls that is both avoidable
and for which $\sum \big(r_n/\vert c_n\vert\big)^{d-2}$ is
divergent: in fact, the possibility that condition \eqref{1} is
redundant in Theorem~\ref{t2} has not been ruled out as yet. We
address this gap in our final result.

\begin{thm}\label{t4}Suppose that $f$ is any increasing
unbounded function on $[0,\infty)$. Then there is a separated and
avoidable collection of balls $\overline{B}(c_n,r_n)$, $n \geq 1$,
for which
\[
r_n^{d-2}\vert c_n\vert^{2} \leq f(\vert c_n\vert)\ \ \mbox{ and }\
\  \sum_{n=1}^\infty \left(\frac{r_n}{\vert c_n\vert}\right)^{d-2} =
\infty.
\]
\end{thm}

We will write $\omega(x,E;D)$ to denote the harmonic measure at $x$
of a subset $E$ on the boundary of a region $D$ with respect to $D$.

\section{Proof of Proposition~\ref{t1}}
\noindent We suppose that \eqref{0} holds and choose $N$ so
large that
\[
\sum_{n=N+1}^\infty \left( \frac{r_n}{\vert c_n\vert} \right)^{d-2}
< \frac{1}{2}.
\]
We write $\Omega_N$ for $\R^d\setminus \bigcup_{n>N}
\overline{B}(c_n,r_n)$. For $n>N$, the harmonic measure of the
sphere $S(c_n,r_n)$ at 0 with respect to $\Omega_N$ is less than its
harmonic measure with respect to the larger domain $\R^d\setminus
\overline{B}(c_n,r_n)$, which is $(r_n/\vert c_n\vert)^{d-2}$. Thus
the combined harmonic measure at $0$ with respect to $\Omega_N$ of
the spheres $S(c_n,r_n)$, $n>N$, is at most $1/2$. As a consequence,
Brownian motion in $\R^d$ starting from 0 has a positive
probability (at least $1/2$) of avoiding the set $E = \bigcup_{n>N}
S(c_n,r_n)$ indefinitely.

We write $u(x)$ for the harmonic measure $\omega(x,E;\Omega_N)$, so
that $u(0) < 1/2$. The set of points $x$ in $\Omega_N$ at which
$u(x)<1/2$ is unbounded. In fact, suppose that it was the case that
$u\geq 1/2$ on $S(0,R)\cap \Omega_N$. We could then apply the
maximum principle to the harmonic function $u$ in
$B(0,R)\cap\Omega_N$, noting that $u=1$ on the boundary of
$\Omega_N$ inside $B(0,R)$ (that is, on $E\cap B(0,R)$), and deduce
that $u\geq 1/2$ in $B(0,R)\cap\Omega_N$.

We now write $F$ for the bounded set $\bigcup_{n\leq N} S(c_n,r_n)$,
and choose $R$ so that $F \subset B(0,R)$. Then, for $\vert x \vert
> R$,
\[
\omega(x,F;\Omega)\ \leq\
\omega\Big(x,S(0,R);\R^d\setminus\overline{B}(0,R)\Big)\ =\ \left(
\frac{R}{\vert x \vert}\right)^{d-2}.
\]
It follows that, as $\vert x\vert \to \infty$ in $\Omega$, the
harmonic measure $\omega(x,F;\Omega)$ tends to 0. Thus we may be
sure that there is a point $x_0$ in $\Omega$ for which both
$\omega(x_0,F;\Omega) < 1/2$ and $\omega(x_0,E;\Omega) \leq
\omega(x_0,E;\Omega_N)< 1/2$. The finite boundary of $\Omega$, that
is $E\cup F = \bigcup_{n=1}^\infty S(c_n,r_n)$,  does not have full
harmonic measure at $x_0$. By the maximum principle, it does not
have full harmonic measure at 0 either, and so the balls
$\overline{B}(c_n, r_n)$, $n\geq 1$, do not hide infinity from the
origin.

\section{Proof of Theorem~\ref{t2}}
\noindent Let us suppose that \eqref{1} holds and that
$\sum_{n=1}^\infty \big(r_n/\vert c_n\vert\big)^{d-2}$ is divergent.
We wish to show that Brownian motion starting from $0$ will never
escape to infinity in $\Omega$.

We set $I_m =  \{n \in \mathbb N\colon \epsilon 2^{m-1}< \vert
c_n\vert \leq \epsilon 2^m\}$. We set
\[
m_0 = \left\lceil \frac{3d-5}{d-2}\right\rceil
\]
and note that there is a $k_0$ between 1 and $m_0$ inclusive for
which
\[
\sum_{j=0}^\infty \sum_{n\in I_{k_0+j m_0}} \left(\frac{r_n}{\vert
c_n\vert}\right)^{d-2} = \infty.
\]
We ignore all balls whose index does not lie in $I_{k_0+j m_0}$ for
some $j$: with fewer balls to avoid, it is easier for Brownian
motion starting from 0 to escape to infinity in this new domain
$\Omega$. The balls that remain lie more or less in annuli whose
inner radius is half the outer radius but arranged so that the
annuli are far apart, in that the inner radius of each annulus is
$2^{m_0}$ times that of the previous annulus.

Following the argument of Ortega-Cerd\`{a} and Seip
\cite[p.~909]{OCS}, we write $m_j$ for $k_0+j m_0$, $R_j$ for
$\epsilon 2^{m_j-1}$, $S_j$ for $S(0,R_j)$ and set $P_j$ to be the
probability that Brownian motion in $\Omega$ starting from $0$ hits
$S_j\cap\Omega$. We need to show that $P_j \to 0$ as $j\to\infty$.

We let $Q_j$ be the supremum of the probabilities that Brownian
motion with starting point on $S_j\cap\Omega$ hits
$S_{j+1}\cap\Omega$. Then
\[
P_{j+1} \leq Q_j \,P_j
\]
and so
\[
P_{n+1} \leq P_1\,\prod_{j=1}^n Q_j.
\]
If $0<a_j <1$ and $\sum_{j=1}^\infty (1- a_j)$ is divergent, then
the infinite product $\prod_{j=1}^\infty a_j = 0$. Theorem~\ref{t2}
therefore follows from the next lemma.

\begin{lemma}\label{l1}
We set $C$ to be $1+4^{d+3}M\epsilon^{-d}$. Then, for all
sufficiently large $j$,
\begin{equation}\label{l1eqn}
1 - Q_j\ \geq\ \frac{1}{2^{d-1}\,C} \sum_{n \in I_{m_j}}\left(
\frac{r_n}{\vert c_n\vert} \right)^{d-2}
\end{equation}
\end{lemma}

\begin{proof}
We write
\[
\Omega_j\ =\ B(0,R_{j+1}) \setminus \bigcup_{n\in I_{m_j}}
\overline{B}(c_n,r_n) .
\]
Then $Q_j \leq \hat{Q}_j$ where $\hat{Q}_j$ is the supremum of the
probabilities that a Brownian motion in $\Omega_j$ with starting
point on $S_j\cap\Omega_j$ hits $S_{j+1}\cap\Omega_j$.
Lemma~\ref{l1} may be proved, therefore, by showing that
\begin{equation}
\inf_{x \in S_j\cap \Omega_j} \omega\Big(x, \bigcup_{n\in
I_{m_j}}S(c_n,r_n)\,;\Omega_j\Big)\ \geq\ \frac{1}{2^{d-1}\,C}
\sum_{n\in I_{m_j}}\left(\frac{r_n}{\vert c_n\vert}\right)^{d-2}
\label{star}
\end{equation}

We consider
\[
u(x)\ =\ \sum_{n\in I_{m_j}} \left( \frac{r_n}{\vert
x-c_n\vert}\right)^{d-2} \  \mbox{ for } x\in \Omega_j,
\]
so that $u$ is harmonic in $\Omega_j$. We suppose that $x\in
S(c_m,r_m)$ for some $m\in I_{m_j}$. Then $r_m/\vert x-c_m\vert =
1$. We now show that the assumption that $r_n^{d-2} \vert c_n\vert^2
\leq M$, for each $n$, leads to
\begin{equation}\label{3.1}
\sum_{\begin{subarray}{c}n\neq m\\ n\in I_{m_j}\end{subarray}}
\frac{r_n^{d-2}}{\vert x-c_n\vert^{d-2}}\ \leq\
4^{d+3}M\epsilon^{-d}.
\end{equation}
By \eqref{1}, we may assume that $r_n < \epsilon$ for $n \in
I_{m_j}$, once $j$ is sufficiently large. The separation condition
\eqref{separation} implies that there are at most $4^d 2^{kd}$ balls whose
centres
lie at a distance of more than $2^k\epsilon$ but less than
$2^{k+1}\epsilon$ from $x$, for $k\geq 1$. Each putative ball in
this annulus contributes at most
\[
\frac{M}{\vert c_n\vert^2} \frac{1}{\vert x-c_n\vert^{d-2}}\ \leq\
\frac{M}{R_j^2}\, \frac{1}{(2^{k}\epsilon)^{d-2}}\ =\
\frac{M4^{k}}{\epsilon^{d-2}}\, \frac{1}{R_j^2\, 2^{kd}}
\]
to the sum in \eqref{3.1}. Since $m_j+1$ annuli centred at $x$ will
cover all balls $B(c_n,r_n)$ with $n \in I_{m_j}$, we find that
\begin{eqnarray*}
\sum_{\begin{subarray}{c} n\neq m\\ n\in I_{m_j}\end{subarray}}
\frac{r_n^{d-2}}{\vert x-c_n\vert^{d-2}}& \leq & \sum_{k=1}^{m_j+1}
4^d\, 2^{kd}\, \frac{M4^k}{\epsilon^{d-2}}\, \frac{1}{R_j^2\,
2^{kd}}\cr
& = &  \frac{4^d\,M}{R_j^2\,\epsilon^{d-2}}
\sum_{k=1}^{m_j+1} 4^k\ \leq\
\frac{4^{d+2}M}{R_j^2\,\epsilon^{d-2}}\,4^{m_j}.
\end{eqnarray*}
Since $R_j^2 = \epsilon^2 4^{m_j-1}$, the estimate \eqref{3.1}
follows. We have shown that the harmonic function $u$ satisfies
\begin{equation}\label{3.2}
u(x) \leq 1+4^{d+3}M\epsilon^{-d} = C\ \mbox{ for }\ x \in
\bigcup_{n\in I_{m_j}}S(c_n,r_n).
\end{equation}
We now need an estimate for the size of $u$ on the sphere $S_{j+1}$.
If $\vert c_n \vert \leq 2R_j$ and $\vert x\vert = R_{j+1}$, then
\[
\vert x-c_n\vert \geq R_{j+1}-2R_j = (2^{m_0}-2)R_j \geq
2^{m_0-1}R_j.
\]
Thus, for $x \in S_{j+1}$,
\begin{eqnarray}
u(x) & = &\sum_{n\in I_{m_j}} \left( \frac{r_n}{\vert
    x-c_n\vert}\right)^{d-2}\nonumber\\
&\leq & \frac{1}{2^{(m_0-1)(d-2)}}\sum_{n\in I_{m_j}}
    \left( \frac{r_n}{R_j}\right)^{d-2}\nonumber\\
& \leq & \frac{4^{d-2}}{2^{m_0(d-2)}}\sum_{n\in I_{m_j}} \left( \frac{r_n}{\vert
    c_n\vert}\right)^{d-2}\nonumber\\
& \leq &  \frac{1}{2^{d-1}}
    \sum_{n\in I_{m_j}} \left( \frac{r_n}{\vert
    c_n\vert}\right)^{d-2}
\label{3.3}
\end{eqnarray}
It follows from \eqref{3.2}, \eqref{3.3} and the maximum principle
that, for $x \in \Omega_j$,
\[
C\,\omega\Big(x, \bigcup_{n\in I_{m_j}}S(c_n,r_n);\Omega_j\Big)\
\geq\ u(x) - \frac{1}{2^{d-1}}
    \sum_{n\in I_{m_j}} \left( \frac{r_n}{\vert
    c_n\vert}\right)^{d-2}
\]
Finally, we use this inequality with $x\in S_j\cap \Omega_j$. For
such $x$, we have $\vert x-c_n\vert \leq \vert x \vert + \vert
c_n\vert \leq 2\vert c_n\vert$, and so
\[
u(x) \geq \frac{1}{2^{d-2}}\sum_{n\in I_{m_j}} \left(
\frac{r_n}{\vert
    c_n\vert}\right)^{d-2}
\]
The estimate \eqref{star} follows immediately. This completes the
proof of the lemma, and hence of Theorem~\ref{t2}.
\end{proof}
\begin{remark}
If the centres of the balls lie on a $(d-1)$-dimensional hyperplane,
then the conclusion of Theorem~\ref{t2} still holds with the
assumption \eqref{1} replaced by the weaker assumption $r_n^{d-2}
\vert c_n\vert \leq M$. Working through the proof of Lemma~\ref{l1},
it is still possible to conclude that $u$ is bounded on the boundary
of the balls with index in $I_{m_j}$ by a constant that is
independent of $j$, as in \eqref{3.2}. (In fact, there are at most
$4^{d-1}2^{k(d-1)}$ balls \lq whose centres lie at a distance of
more than $2^k\epsilon$ but less than $2^{k+1}\epsilon$ from $x$,
for $k\geq 1$\rq.) The remainder of the proof of Lemma~\ref{l1} is
unchanged.

If the centres of the balls lie on a $(d-2)$-dimensional hyperplane
then one may replace \eqref{1} by the weaker assumption
$r_n^{d-2}\log \vert c_n\vert \leq M$ and still retain the
conclusion of Theorem~\ref{t2}. For example, suppose that in $\R^3$
we put a ball of radius $r_n$ at the point $(n,0,0)$, for $n \geq
2$. Under the assumption that $r_n \leq M/\log n$, this string of
beads in $\R^3$ is avoidable if and only if $\sum r_n/n$ is finite.

If the centres of the balls lie on a $(d-3)$-dimensional hyperplane,
then it suffices to assume the the radii of the balls are uniformly
bounded in order for Theorem~\ref{t2} to hold. In this case,
however, there can be at most about $m^{d-4}$ balls whose distance
from the origin is about $m$. Assuming that the radii of the balls
are bounded by $R$, say, it follows that
\[
\sum_{n}\left( \frac{r_n}{\vert c_n\vert}\right)^{d-2} \leq R^{d-2}
\sum_n \frac{1}{\vert c_n\vert^{d-2}} \leq  C R^{d-2} \sum_m
m^{d-4}\frac{1}{m^{d-2}}
\]
which is finite. A collection of balls of uniformly bounded radius
whose centres lie on a $(d-3)$-dimensional hyperplane will always be
avoidable.
\end{remark}
\section{Proof of Theorem~\ref{t3}}
\noindent

To begin with we note that, in the case of a regularly located
configuration of balls, the sum $\sum_n\left(r_n/\vert
c_n\vert\right)^{d-2}$ and the integral $ \int^\infty r
\phi(r)^{d-2}\,dr$ are comparable. The implication that the balls
are avoidable if $ \int^\infty r \phi(r)^{d-2}\,dr$ is finite is now
an immediate consequence of Proposition~\ref{t1}.

The reverse implication will follow from Theorem~\ref{t2} once we
check that the condition \eqref{1} is automatically satisfied under
the regularity assumption if the balls are avoidable. We establish
this in the next lemma, whose proof bears a certain resemblance to
that of Lemma~\ref{l1}.

\begin{lemma}\label{l2}
Suppose that the balls $\{\overline{B}(c_n,r_n))\}_{n\geq 1}$ are
regularly located and that $r^2 \phi(r)^{d-2}$ is an unbounded
function of $r$. Then the collection of balls is not avoidable.
\end{lemma}

\begin{proof}
There is a sequence of radii $\{R_j\}_{j=1}^\infty$ for which $R_j^2
\phi(2R_j)^{d-2} \to \infty$ as $n \to \infty$. We put
\[
C = \frac{A_2}{2A_3}
\]
where the particular numbers $A_2$ and $A_3$ that we need depend on
the dimension, on the separation number $\epsilon$ and on the
density number $R$ but on nothing else, and may be worked out in
principle from the proof that follows. We assume that $R_{j+1}>4R_j$
and that $(R_j/R_{j+1})^{d-2} \leq C$ for each $j$. For a technical
reason, we change the definition of $\phi$ in the following way: we
set $\widetilde\phi(x)=\phi(2R_j)$ if $x\in[R_j,2R_j]$ for some $j$
and $\widetilde \phi(x)=\phi(x)$ elsewhere. We take new balls
$B(c_n,\widetilde\phi(|c_n|))$. The size of the balls is thereby
decreased: thus if the new balls are unavoidable then the original
ones are unavoidable too. For the sake of simplicity, we will still
denote by $\phi$ the regularized $\widetilde \phi$ and the new
smaller balls will still be called $B(c_n,r_n)$. We write $S_j$ for
the sphere $S(0,R_j)$ and $\phi_j = \phi(R_j) = \phi(2R_j)$.

Arguing as in the proof of Theorem~\ref{t2}, we let $Q_j$ be the
supremum of the probabilities that Brownian motion in $\Omega$ with
starting point on $S_j\cap\Omega$ hits $S_{j+1}\cap\Omega$, and wish
to show that $\prod_{j=1}^\infty Q_j =0$, that is that
\[
\sum^\infty_{j=1} (1-Q_j) = \infty.
\]
We write  $I_j$ for $\big\{n: R_j \leq \vert c_n\vert \leq
2R_j\big\}$, and write
\[
\Omega_j\ =\ B(0,R_{j+1}) \setminus \bigcup_{n\in I_j}
\overline{B}(c_n,r_n).
\]
Then $Q_j$ is bounded above by $\hat{Q}_j$, the supremum of the
probabilities that Brownian motion with starting point on
$S_j\cap\Omega_j$ hits $S_{j+1}\cap\Omega_j$. We will show that, for
all sufficiently large $j$,
\begin{equation}\label{4.00}
1-\hat{Q}_j = \inf_{x \in S_j\cap\Omega_j} \omega\Big(x,
\bigcup_{n\in I_j}S(c_n,r_n);\Omega_j\Big)\ \geq\ \delta,
\end{equation}
for some positive $\delta$. Again we consider
\[
u(x)\ =\ \sum_{n\in I_j} \left( \frac{r_n}{\vert
x-c_n\vert}\right)^{d-2}, \ \ x\in \Omega_j,
\]
so that $u$ is harmonic in $\Omega_j$. Since $\phi$ is constant on
$[R_j, 2R_j]$, we have $r_n  = \phi(R_j) = \phi_j$ for $n \in I_j$
and
\[
u(x)\ =\ \phi_j^{d-2} \sum_{n\in I_j} \frac{1}{\vert
x-c_n\vert^{d-2}}.
\]
Suppose that $x$ lies on the boundary of a ball $S(c_m,r_m)$ with $m
\in I_j$. It is a consequence of the separation condition that there
can be at most $A2^{kd}$ balls with centre at a distance that is
between $\epsilon 2^{k-1}$ and $\epsilon 2^k$ from $x$, with $k\geq
1$. Each such ball contributes at most $A2^{-k(d-2)}$ to the above
sum, making for a combined contribution of at most $A2^{2k}$. We
need only count those $k$ with $\epsilon 2^k \leq 6R_j$, as there
are no balls under consideration that are more distant than $6R_j$
from $x$. The ball $B(c_m,r_m)$ itself contributes 1 to $u(x)$,
which leads to the estimate
\[
u(x)\leq  1+\phi_j^{d-2} \sum_{k\colon\epsilon2^k\leq 6R_j} A2^{2k}
\leq 1+ AR_j^2 \phi_j^{d-2}.
\]
As $R_j^2 \phi_j^{d-2} \geq 1$ for sufficiently large $j$,
\begin{equation}\label{4.1}
u(x)\leq A_1R_j^2 \phi_j^{d-2}\ \mbox{ for }\ x \in \bigcup_{n\in
I_j} S(c_n,r_n).
\end{equation}
Here $A_1$ is some appropriate number that depends only on the
dimension and on the separation number $\epsilon$.

For $x\in S_j$, we have $\vert x-c_n\vert \leq 4R_j$. At this point
we use the assumption that the balls are uniformly dense to deduce
that
\begin{eqnarray*}
u(x)&= &\phi_j^{d-2} \sum_{n\in I_j} \frac{1}{\vert
x-c_n\vert^{d-2}} \nonumber \\
&\geq& \phi_j^{d-2} \frac{1}{(4R_j)^{d-2}}\,\sum_{n\in I_j} 1\nonumber \\
&\geq& \phi_j^{d-2} \frac{A_2}{R_j^{d-2}}\,R_j^d
\end{eqnarray*}
where the number $A_2$ depends only on the dimension and on the
number $R$ that appears in the definition of \lq regularly
located\rq. Thus,
\begin{equation}\label{4.2}
u(x) \geq A_2 R_j^2\,\phi_j^{d-2} \ \mbox{ for }\ x \in S_j.
\end{equation}

Finally, for $x$ on the sphere $S_{j+1}$ and $n \in I_j$, we have
$\vert x-c_n\vert \geq R_{j+1} - 2R_j \geq R_{j+1}/2$. Hence on
$S_{j+1}$ the function $u$ satisfies
\[
u(x)\ \leq\ \phi_j^{d-2} \frac{2^{d-2}}{R_{j+1}^{d-2}}\,\sum_{n\in
I_j} 1\ \leq\ A_3 \phi_j^{d-2} \frac{R_j^d}{R_{j+1}^{d-2}}.
\]
Since $(R_j/R_{j+1})^{d-2} \leq C$, we obtain that
\begin{equation}\label{4.3}
u(x)  \leq \frac{1}{2}A_2\, R_j^2 \phi_j^{d-2}\ \mbox{ for }\ x \in
S_{j+1}.
\end{equation}
It follows from \eqref{4.1}, \eqref{4.3} and the maximum principle
that, for $x \in \Omega_j$,
\[
A_1 R_j^2 \phi_j^{d-2} \,\omega\Big(x, \bigcup_{n\in
I_j}S(c_n,r_n);\Omega_j\Big)\ \geq\ u(x) - \frac{1}{2} A_2\, R_j^2
\phi_j^{d-2}
\]
Making use of \eqref{4.2}, we deduce from this last estimate that,
for $x \in S_j$,
\[
A_1 \,\omega\Big(x, \bigcup_{n\in I_j}S(c_n,r_n);\Omega_j\Big)\
 \geq\ \frac{1}{2}A_2.
\]
Thus \eqref{4.00} has been proven. \end{proof}

\begin{remark}
With the same proof, one may consider a slightly more general
situation where one changes the metric. Assume that a function
$\psi:\R^+\to\R^+$ satisfies the smoothness condition $\psi(y)\simeq
\psi(x)$ whenever $ y< x<2y$. We say that a sequence $\{c_n\}$  is
$\psi$-regularly located if there is a $\delta>0$ such that the
balls $B(c_n,\delta\psi(|c_n|))$ are pairwise disjoint and there is
an $R>0$ such that any ball $B(x,R\psi(x))$ contains at least a
center $c_n$. Assume finally that we have a sequence of disjoint
balls with $\psi$-regularly located centres and the radii of the
balls depend on the centre, $r_n=\phi(|c_n|)$, where $\phi$ is a
decreasing positive function. Then the balls are avoidable if and
only if
\[
\int^{+\infty} \frac{x\phi(x)^{d-2}}{\psi(x)^{d}}\, dx <\infty.
\]
The case $\psi=1$ is the case previously considered.
\end{remark}

\section{Construction of the examples: Proof of Theorem~\ref{t4}}
\noindent

We wish to show by examples that the assumption \eqref{1} in
Theorem~\ref{t2} is necessary. The examples are of avoidable and
separated configurations of balls for which the series
$\sum\big(r_n/\vert c_n\vert\big)^{d-2}$ is divergent, in which case
$r_n^{d-2}\vert c_n\vert^2$ must be unbounded by Theorem~\ref{t2}.
In Theorem~\ref{t4} it is asserted that such configurations of balls
are possible even with a growth restriction on $r_n^{d-2}\vert
c_n\vert^2$. Leaving the growth restriction to one side for the
moment, we first give the details of a plain vanilla example that
incorporates the idea behind the general construction.

\begin{proposition}
There is an avoidable, separated configuration of balls,
$\overline{B}(c_n,r_n)$, $n\geq 1$, in $\R^3$ for which
\[
\sum_{n=1}^\infty \frac{r_n}{\vert c_n\vert} = \infty
\]
\end{proposition}

\begin{proof}
Consider a string of closed balls $\overline{B}_1$,
$\overline{B}_2$, $\ldots$ $\overline{B}_{2k}$, each of radius $1/4$
and with centres $c_i$ on the $x_1$-axis at $m+i$, $i=1$, $2$,
$\ldots$ $2k$. We write $S(m,k) = \cup_{i=1}^{2k} \overline{B}_i$
and wish to estimate $\omega\big(0,\partial S(m,k);\R^3\setminus
S(m,k)\big)$. We consider, as ever,
\[
u(x) = \sum_{i=1}^{2k} \frac{1}{\vert x-c_i\vert}
\]
Suppose that $x$ lies on the boundary of one of the balls
$\overline{B}_j$. Then $\vert x-c_i\vert \leq \vert i-j\vert +1/4
\leq 2\vert i-j\vert$ for $i\ne j$, and there are at least $k$ balls
to one side or other of any one ball. It follows that
\[
u(x) \geq \frac{1}{2} \sum_{i=1}^k \frac{1}{i} \geq \frac{1}{2} \log
k, \quad x \in \partial S(m,k).
\]
By the maximum principle,
\[
\omega\big(0,\partial S(m,k);\R^3\setminus S(m,k)\big) \leq
\frac{2}{\log k}\, u(0) = \frac{2}{\log k} \sum_{i=1}^{2k}
\frac{1}{m+i} \leq \frac{4k}{m\log k}.
\]
We construct our counterexample as follows. Let $S_n = S(n^2,
 \lfloor n/\log n \rfloor)$ and
\[
\Omega = \R^3 \setminus \bigcup_{n=n_0}^\infty S_n.
\]
Then
\begin{eqnarray*}
\omega(0,\partial \Omega,\Omega) & \leq & \sum_{n=n_0}^\infty
\omega\big(0,\partial S_n;\R^3\setminus S_n\big)\cr
& \leq & 4\sum_{n=n_0}^\infty \frac{n/\log n}{n^2 \log\lfloor n/\log n
\rfloor}\cr
& \leq & 48\sum_{n=n_0}^\infty \frac{1}{n\log^2 n}
\end{eqnarray*}
Thus $n_0$ may be chosen to be sufficiently large so that the balls
are separated and so that $\omega(0,\partial \Omega,\Omega) <1$, in
which case the balls are avoidable.

On the other hand, the contribution of each string of balls $S_n$ to
the series $\sum r_n/\vert c_n\vert$ is comparable to $1/(n \log
n)$, and this sum is divergent. \end{proof}

In the examples that follow the balls are arranged in clusters
rather than in higher dimensional strings, though the reason the
examples work is the same: each ball in a cluster of balls
contributes significantly less to the harmonic measure than it would
do if taken individually.

\begin{proof}[Proof of Theorem~\ref{t4}]We consider a cluster of
$k^d$ balls, each of radius $r$ less than $1/4$, whose centres have
integer coordinates and are evenly distributed in a large ball that
has radius approximately $k$ and is centred at a distance $m$ from
the origin. We assume that $k \leq m/2$ and refer to this cluster of
balls as $C(m,k,r)$. We again use the function
\begin{equation}\label{ucluster}
u(x) = \sum_i  \frac{1}{\vert x-c_i\vert^{d-2}},
\end{equation}
the $c_i$ being the centres of the balls. If $x$ is a point on the
boundary of one of these balls and $1\leq i \leq k$, there are at
least $a\,i^{d-1}$ balls whose centres lie at a distance at most
$2i$ from $x$. Here $a$ represents a number that depends only on the
dimension. Moreover, no ball needs to be chosen twice, that is for
two different values of $i$. We find that, for a point $x$ on the
boundary of any ball in the cluster,
\[
u(x) \geq \sum_{i=1}^k \frac{1}{(2i)^{d-2}} \,a i^{d-1} \geq a k^2.
\]
By the maximum principle,
\[
\omega\big(0,\partial C(m,k,r),\R^d\setminus C(m,k,r)\big) \leq
\frac{A}{k^2}\,u(0) \leq A\,\frac{k^d}{m^{d-2}}\,\frac{1}{k^2} =
A\left( \frac{k}{m}\right)^{d-2}
\]

We suppose that an increasing unbounded function $f$ on $[0,\infty)$
is given. To each positive integer $n$ there corresponds a choice of
variable $m_n$ for which $f(m_n) \geq n^{2d}$ and $m_n > 2m_{n-1}$.
We then choose $k_n$ to be $m_n/n^2$ and choose the radius $r_n$ so
that $r_n^{d-2}m_n^2 = f(m_n)$. [We assume that the function $f$
satisfies $f(x) \leq 4^{2-d}x^2$, so that $r_n <1/4$.] We set
\[
\Omega = \R^d \setminus \bigcup_{n=n_0}^\infty C(m_n,k_n,r_n)
\]
and write $\omega_n(x)$ for the harmonic measure at $x$ of the
finite boundary of $C(m_n,k_n,r_n)$ with respect to the domain $\R^d
\setminus C(m_n,k_n,r_n)$. Then
\[
\omega(0,\partial \Omega,\Omega) \leq \sum_{n=n_0}^\infty w_n(0)
\leq \sum_{n=n_0}^\infty A\left( \frac{k_n}{m_n}\right)^{d-2}= A
\sum_{n=n_0}^\infty \frac{1}{n^{2(d-2)}},
\]
which we can arrange to be strictly less than 1 by taking $n_0$ to
be sufficiently large.

The sum in \eqref{0} for this collection of balls is comparable to
\[
\sum_{n=n_0}^\infty k_n^d \left(\frac{r_n}{m_n}\right)^{d-2}
\]
Since $r_n^{d-2}m_n^2 = f(m_n) \geq n^{2d}$, the general term in
this last sum exceeds $n^{2d} (k_n/m_n)^d$, which in turn exceeds
$n^{2d}(1/n^2)^d = 1$. The sum in \eqref{0} is therefore divergent.
\end{proof}

\section{Addendum: the union of two avoidable sets is avoidable}
\noindent

At a certain point in our research, it seemed that it might be
helpful to know if the union of two avoidable collections of balls
would again be avoidable. Put another way, is it possible to split
an unavoidable collection of balls into two disjoint avoidable
collections? Though the solution to this problem is no longer an
essential ingredient in the proofs we have presented here, we cannot
resist including the elegant solution to this problem found by
Professor Rosay. We are grateful to him for granting us permission
to include his proof in this article.

A set $A$  is called \emph{avoidable from} $p$ if Brownian motion in
$\R^d$ starting at $p$ has a probability smaller than one of hitting
$A$. We assume that $\R^d\setminus A$ is connected: then, by the
maximum principle, if $A$ is avoidable from one point it is
avoidable from any other point. In this case we just say that the
set $A$ is \emph{avoidable}. Equivalently $A$ is avoidable whenever
there is a positive  harmonic function $u$ in $\R^d\setminus A$ such
that $u\equiv 1$ on the boundary of $A$ but $\inf u = 0$.
\begin{proposition}\label{p3}
If two avoidable sets $A$ and $B$ satisfy\/  $\R^d\setminus(A\cup
B)$ is connected then $A\cup B$ is avoidable.
\end{proposition}
\noindent The basic lemma required to prove this proposition is the
following:
\begin{lemma}
If $E$ is avoidable and $u_E$ is the associated positive harmonic
function in $\R^d\setminus E$, with $u_E\equiv 1$ on the boundary of
$E$ and  $\inf u_E = 0$, then there is an $R_0$ such that for all
$R\ge R_0$ the set of points
\[
S^R_E=\{ x\in S(0,R)\setminus E\, \colon\,u_E(x)\le 1/4\}
\]
satisfies $|S_E^R|> \frac{3}{4}|S(0,R)|$. Here the measure indicated
by $\vert \cdot \vert$ is Lebesgue area measure on $S(0,R)$.
\end{lemma}
\begin{proof}
We take a point $q$ where $u_E(q)<1/32$. For any  $R$ with $R>|q|$,
we denote by $\mu_R$ the harmonic measure on the boundary of
$B(0,R)\setminus E$ with respect to $q$. Then, since $u = 1$ on
$\partial E\cap B(0,R)$, we have
\[
\frac{1}{32} >  \mu_R\big(\partial E\cap B(0,R)\big) +
\frac{1}{4}\big( 1-\mu_R\big(\partial E\cap B(0,R)\big) - \mu_R
(S^R_E)\big)
\]
from which it follows that $\mu_R (S^R_E) > 7/8$. We denote by
$\sigma_R$ harmonic measure with base point $q$ with respect to the
ball $B(0,R)$, so that $\sigma_R\ge \mu_R$ on $S(0,R)$. Thus,
$\sigma_R(S^R_E) > 7/8$ for all $R>|q|$. The harmonic measure
$\sigma_R$ can be given explicitly, but the key property is that as
$R\to\infty$ it is more and more similar to the normalized area
measure on $S(0,R)$. Thus $|S_E^R|>\frac{3}{4}|S(0,R)|$ for all
large $R$.
\end{proof}
\begin{proof}[Proof of Proposition~\ref{p3}] For the sets $A$ and $B$ we take
the
corresponding functions $u_A$ and $u_B$. We take $R$ so that
$|S_A^R|>\frac{3}{4}|S(0,R)|$ and $|S_B^R|>\frac{3}{4}|S(0,R)|$.
This means that there is point $p$ that lies in the intersection
$S^R_A\cap S^R_B$. We define $u=u_A+u_B$: it is a positive and
bounded harmonic function defined outside $A\cup B$. On the boundary
of $A\cup B$ it satisfies $u\ge 1$ and on the other hand $u(p) \leq
1/2$. Thus $A\cup B$ is avoidable from $p$. Since the complement
$\R^d\setminus(A\cup B)$ is connected, then it is avoidable from any
point.
\end{proof}

\end{document}